\documentclass{article}

\title{Small ball probability estimates\\ in terms of width
\thanks{Research partially supported by KBN Grant 2 PO3A 027 22}}
\author{Rafa{\l} Lata{\l}a and Krzysztof Oleszkiewicz}
\date{}

\newtheorem{Th}{Theorem}
\newtheorem{Cor}{Corollary}
\newtheorem{Conj}{Conjecture}
\newtheorem{lem}{Lemma}

\def\er{{\mathbb R}}

\def\Pr{{\mathbf P}}
\def\Ex{{\mathbf E}}

\usepackage{amssymb}

\begin{document}
\maketitle

\begin{abstract}
A certain inequality conjectured by Vershynin is studied.
It is proved that for any symmetric convex body $K \subseteq \er^{n}$
with inradius $w$ and $\gamma_{n}(K) \leq 1/2$ there is
\[
\gamma_{n}(sK) \leq (2s)^{w^{2}/4}\gamma_{n}(K)
\]
for any $s \in [0,1].$ Some natural corollaries are deduced. Another
conjecture of Vershynin is proved to be false.
\end{abstract}

\section{Introduction}
In his lecture at Snowbird'2004 AMS Conference Roman Vershynin posed
two conjectures related to the rate of decay of the Gaussian measure of convex
symmetric sets under homothetic shrinking. The first conjecture concerned
some bounds in terms of the width of a convex symmetric set. The second one
stated that among all convex symmetric bodies with fixed both Gaussian measure
and width the cylinders are the ones that have the slowest decay of Gaussian measure under
homothetic shrinking. Both conjectures will be described more precisely below.
In this paper we prove some version of the first conjecture (Section 2)
and we demonstrate that the second conjecture cannot hold in general
(Section 4). We also sketch some natural direct applications that motivated
Vershynin's questions. More sophisticated geometric consequences related
to the Dvoretzky theorem were recently proved by Klartag and Vershynin,
\cite{KV}.

Let us introduce some notation and results that will be used.
By $\gamma_{n}$ we will denote
the standard Gaussian probability measure on $\er^{n},$
with $\gamma_{n}(dx)=(2\pi)^{-n/2}e^{-|x|^{2}/2}dx.$
For a set $A$ in $\er^{n}$ we will write
$A_{t}=\{x\in \er^{n}\colon d(x,A)<t\}$.
The Gaussian isoperimetry \cite{B,ST} states that if
$\gamma_{n}(A)=\Phi(x)$, then
$\gamma_{n}(A_{t})\geq \Phi(x+t)$ for all $t>0$.

The S-inequality \cite{LO} says that if $K\subset\er^{n}$ is convex
symmetric with
$\gamma_{n}(K)=\gamma_{1}([-a,a])$ then
$\gamma_{n}(tK)\leq \gamma_{1}([-ta,ta])$ for all $0<t<1$. This easily implies
that for some universal constant $C$
\begin{equation}
\label{estdil1}
\forall_{0\leq t\leq 1}\ \gamma_{n}(tK)\leq Ct\gamma_{n}(K) \mbox{ if }
\gamma_{n}(K)\leq \frac{1}{2}.
\end{equation}
The inequality (\ref{estdil1}) was first proved in \cite{6A}
(see also its generalization to log-concave measures in \cite{L}).
Although in general one cannot improve (\ref{estdil1}),
Vershynin conjectured that a stronger inequality can hold
for the sets of large width.

For a convex symmetric set $K$ in $\er^{n}$ let the {\em inradius of $K$} be defined as
\[w(K)=\sup\{r>0\colon B(0,r)\subset K\}.\]
Notice also that $w(K)$ is half of the {\em width of $K$}.

The B-inequality proved recently by Cordero, Fradelizi and Maurey
\cite{CFM} gives that for any symmetric
convex set $K$ in $\er^{n}$ the function $t\mapsto \ln \gamma_{n}(e^{t}K)$ is concave.
In particular
\begin{equation}
\label{estcfm}
 \forall_{0<s\leq t\leq 1}\
 \frac{\gamma_{n}(sK)}{\gamma_{n}(K)}\leq (\frac{\gamma_{n}(tK)}{\gamma_{n}(K)})^{\frac{\ln s}{\ln t}}.
\end{equation}

\section{Main Results.}

\begin{Th}
For any convex symmetric set $K$ in $\er^{n}$ with $\gamma_{n}(K)\leq 1/2$
we have
\begin{equation}
\label{est1}
  \forall_{0\leq t\leq \frac{1}{2}}\ \gamma_{n}(tK)\leq
  t^{\frac{\ln 2}{8}w^{2}(K)}\gamma_{n}(K).
\end{equation}
\end{Th}

{\bf Proof.} Let $s\geq 1$ be such that $\gamma_{n}(sK)=1/2$. By
the concavity of $t\mapsto\gamma_{n}(e^{t}K)$, we get
$\gamma_{n}(tK)/\gamma_{n}(K)\leq \gamma_{n}(tsK)/\gamma_{n}(sK)$.
Since $w(sK)\geq w(K)$, hence
we may and will assume that $\gamma_{n}(K)=1/2$.

Notice that $\frac{1}{2}K+\frac{1}{2}B(0,w(K))\subset \frac{1}{2}K+\frac{1}{2}K=K$, so $(K^{c})_{w(K)/2}\cap \frac{1}{2}K=\emptyset$. Thus by the Gaussian isoperimetry
\[\gamma_{n}(\frac{1}{2}K)\leq 1-\gamma_{n}((K^{c})_{w(K)/2})\leq
1-\Phi(w(K)/2)\leq \frac{1}{2}e^{-\frac{w^{2}(K)}{8}}\]
\[=\gamma_{n}(K)(\frac{1}{2})^{\frac{\ln 2}{8}w^{2}(K)}.\]
We get (\ref{est1}) by (\ref{estcfm}).\hfill$\Box$

\medskip

The following related conjecture seems reasonable.

\begin{Conj}
\label{hipoteza}
For any $\kappa \in (0,1)$ there exist positive constants $C=C(\kappa)$
and $w_{0}=w_{0}(\kappa)$ such that
for any convex symmetric set $K$ in $\er^{n}$ with $\gamma_{n}(K)\leq 1/2$
and $w(K) \geq w_{0}$
we have
\begin{equation}
\label{conjv}
\forall_{0\leq t\leq 1}\ \gamma_{n}(tK)\leq (Ct)^{\kappa w^{2}(K)}\gamma_{n}(K).
\end{equation}
\end{Conj}

The first Vershynin's conjecture (it was originally formulated
in the language of Theorem \ref{procesy}) was that the above inequality
is true for some fixed $\kappa \in (0,1).$
Inequality (\ref{est1}) shows that (\ref{conjv}) holds for
$\kappa=\frac{1}{8}\ln 2>\frac{1}{12}$,
we will now present some more elaborate argument for $\kappa=\frac{1}{4}$.

\begin{lem}
For any $u,v\geq 0$ we have
\[\gamma_{1}((u+v,\infty))\leq e^{-uv} \gamma_{1}((u,\infty)).\]
\end{lem}

{\bf Proof.} We have
\[\int_{u+v}^{\infty}e^{-s^{2}/2}ds=\int_{u}^{\infty}e^{-(s+v)^{2}/2}ds
  \leq \int_{u}^{\infty}e^{-sv}e^{-s^{2}/2}ds\leq
  e^{-uv}\int_{u}^{\infty}e^{-s^{2}/2}ds.\]

\vspace{-1.5ex}\hfill$\Box$

\begin{Th}
\label{estgauss2}
For any convex symmetric set $K$ in $R^{n}$ with $w=w(K)$ and
$\gamma_{n}(K)\leq \frac{1}{2}$ we have
\[\gamma_{n}(sK)\leq (2s)^{w^{2}/4}\gamma_{n}(K)\ \mbox{ for } s\in[0,1].\]
\end{Th}

{\bf Proof.} Let us notice that $\gamma_{n}(K^{c})\geq \frac{1}{2}$ and
$\frac{1}{2}K\cap (K^{c})_{w/2}=\emptyset$, hence by the isoperimetry
\[\gamma_{n}(\frac{1}{2}K)\leq 1-\gamma_{n}((K^{c})_{w/2})\leq
  \gamma_{1}((w/2,\infty)).\]
Let us define $u\geq w/2$ by the formula
\[\gamma_{n}(\frac{1}{2}K)=\gamma_{1}((u,\infty)).\]

We also have $\frac{t}{2}K\cap ((\frac{1}{2}K)^{c})_{\frac{1-t}{2}w}=\emptyset$
for $0<t<1$,
so again by the isoperimetry and Lemma 1
\[\gamma_{n}(\frac{t}{2}K)\leq \gamma_{1}((u+\frac{1-t}{2}w,\infty))
  \leq e^{-\frac{1-t}{2}wu}\gamma_{1}((u,\infty))\leq
   e^{-(1-t)w^{2}/4}\gamma_{n}(\frac{1}{2}K).\]

Thus
\[\frac{\gamma_{n}(\frac{t}{2}K)}{\gamma_{n}(\frac{1}{2}K)}\leq
  (\frac{t/2}{1/2})^{\frac{w^{2}}{4}\frac{t-1}{\ln t}}.\]
Hence, by the B-inequality for any $s\leq t/2$ we obtain
\[\frac{\gamma_{n}(sK)}{\gamma_{n}(\frac{1}{2}K)}\leq
  (2s)^{\frac{w^{2}}{4}\frac{t-1}{\ln t}}.\]
Taking the limit $t\rightarrow 1^{-}$ we get for $s\in [0,1/2]$
\[\gamma_{n}(sK)\leq\gamma_{n}(\frac{1}{2}K)(2s)^{w^{2}/4}.\]

\vspace{-1.5ex}\hfill$\Box$

Before stating the next result let us introduce some notation. By
$\sigma_{n-1}$ we will denote the (normalized) Haar measure on
$S^{n-1}=\{x\in\er^{n}: |x|=1\}$. For a set $A$ in $\er^{n}$ to
simplify the notation we will write $\sigma_{n-1}(A)$ instead
of $\sigma_{n-1}(A\cap S^{n-1})$.

\begin{Th}
\label{dilsph}
For any convex symmetric set $K$ in $\er^{n}$ with $\sigma_{n-1}(K)\leq
1/2$ we have
\begin{equation}
\label{estsph}
  \forall_{0\leq t\leq 1}\ \sigma_{n-1}(tK)\leq
  (12t)^{\frac{1}{4}(\sqrt{n}w(K)-6)_{+}^{2}}.
\end{equation}
\end{Th}

The proof is based on the following simple lemma.

\begin{lem}
There exists a universal constant $\alpha>1/60$ such that for any
star body $K \subseteq \er^{n}$ we have
\begin{equation}
\label{comp1}
 \gamma_{n}(\sqrt{n}K)\geq \alpha\sigma_{n-1}(K)
\end{equation}
and
\begin{equation}
\label{comp2}
  \gamma_{n}((\sqrt{n}K)^{c})\geq \alpha\sigma_{n-1}(K^{c})
\end{equation}
\end{lem}

{\bf Proof.} Let
\[L=\{\xi\in\er^{n}\colon |\xi|\leq \sqrt{n},\xi/|\xi|\in
  K\cap S^{n-1}\}.\]
Notice that $L\subseteq \sqrt{n}K$ and by the rotational invariance
of $\gamma_{n}$
\[\gamma_{n}(\sqrt{n}K)\geq \gamma_{n}(L)=
  \gamma_{n}(B(0,\sqrt{n}))\sigma_{n-1}(K).\]
Let $X=\sum_{i=1}^{n}(g_{i}^{2}-1)$, where $g_{i}$ are iid
${\mathcal N}(0,1)$ r.v.'s. Since $\Ex|X|/2=\Ex X_{-}\leq
(\Ex X^{2})^{1/2}\Pr(X\leq 0)^{1/2}$ we get
\[\gamma_{n}(B(0,\sqrt{n}))=\Pr(X\leq 0)\geq
\frac{(\Ex X)^{2}}{4\Ex X^{2}}\geq
\frac{(\Ex X^{2})^{2}}{4\Ex X^{4}}\geq \frac{1}{60},\]
where the last inequality follows by an easy calculation, since
$\Ex(g_{i}^{2}-1)=0$, $\Ex(g_{i}^{2}-1)^{2}=2$ and
$\Ex(g_{i}^{2}-1)^{4}=60$.

In a similar way we show that for
\[\tilde{L}=\{\xi\in\er^{n}\colon |\xi|\geq \sqrt{n},\xi/|\xi|\notin
  K\cap S^{n-1}\}\]
we have
\[\gamma_{n}((\sqrt{n}K)^{c})\geq \gamma_{n}(\tilde{L})=
 \Pr(X\geq 0)\sigma_{n-1}(K^{c})\geq\frac{1}{60}\sigma_{n-1}(K^{c}).
 \]

\vspace{-1.5ex}\hfill$\Box$

{\bf Proof of Theorem \ref{dilsph}.} Obviously we may assume that
$\sqrt{n}w(K)>6$ and $t<1/12.$ Let $\alpha$ be the constant given
by the preceding lemma.
If $\sigma_{n-1}(K)\leq 1/2$, then
by (\ref{comp2})
$\gamma_{n}((\sqrt{n}K)^{c})\geq \alpha/2 >\Phi(-2.5).$
Hence by the Gaussian isoperimetry (since
$\mathrm{dist}((\sqrt{n}-s)K,(\sqrt{n}K)^{c})\geq sw(K)$) we get
\[\gamma_{n}((\sqrt{n}-\frac{5}{w(K)})K)\leq
1-\gamma_{n}(((\sqrt{n}K)^{c})_{5})\leq 1-
\Phi(\Phi^{-1}(\alpha/2)+5)\]
\[\leq 1-\Phi(2.5) \leq\min(\alpha, 1/2).\]
Let $t_{0}=\sqrt{n}-5(w(K))^{-1}\geq \sqrt{n}/6>t\sqrt{n}$, then by
(\ref{comp1}) and Theorem \ref{estgauss2}
we get
\[\sigma_{n-1}(tK)\leq \alpha^{-1}\gamma_{n}(\sqrt{n}tK)\leq
  \alpha^{-1}(2t\sqrt{n}/t_{0})^{w^{2}(t_{0}K)/4}\gamma_{n}(t_{0}K)\leq
  (12t)^{(\sqrt{n}w(K)-5)_{+}^{2}}.\]

\vspace{-1.5ex}\hfill$\Box$

\section{Some applications}

Let $X$ be a centered Gaussian vector with values in a separable
Banach space $(F,\| \cdot \|)$. We define $M={\mathrm{Med}}(\|X\|)$ and
\[\sigma=\sigma_{X}=\sup_{f\in E^{*},\|f\|\leq 1}(\Ex f^{2}(X))^{1/2}.\]

\begin{Th}
\label{procesy}
For any $t\in [0,1]$ we have
\[\Pr(\|X\|\leq tM)\leq \frac{1}{2}(2t)^{\frac{M^{2}}{4\sigma^{2}}}.\]
\end{Th}

{\bf Proof.} We may represent $X$ (cf. \cite{LT}) as
$X=\sum_{i=1}^{\infty}x_{i}g_{i}$ for some $x_{i}\in F$ and $g_{i}$ iid ${\mathcal N}(0,1)$ r.v.'s. Standard approximation argument shows that it
is enough to consider the case $X=\sum_{i=1}^{n}x_{i}g_{i}$. Let
\[K:=\{\xi\in\er^{n}\colon \|\sum_{i=1}^{n}x_{i}\xi_{i}\|\leq M\}.\]
Obviously $\gamma_{n}(K)= 1/2$, moreover if $\xi\notin K$ then
for some $f\in E^{*}$, $\|f\|\leq 1$, hence we get
\[M<f(\sum_{i=1}^{n}x_{i}\xi_{i})\leq
 (\sum_{i=1}^{n}f(x_{i})^{2})^{1/2}|\xi|\leq \sigma|\xi|.\]
Thus $w(K)\geq M/\sigma$ and  by Theorem \ref{estgauss2}
\[\Pr(\|X\|\leq tM)=\gamma_{n}(tK)\leq \frac{1}{2}(2t)^{\frac{M^{2}}{4\sigma^{2}}}.\]

\vspace{-1.5ex}\hfill$\Box$

Before formulating the next result let us recall that we define for
$p\neq 0$ the $p$--th moment of a random vector $X$ as
$\|X\|_{p}=(\Ex\|X\|^{p})^{1/p}$ and for $p=0$ as
$\|X\|_{0}=\exp(\Ex \ln \|X\|)$.

\begin{Cor}
For any $p>q>\min(-1,-\frac{M^{2}}{4\sigma^{2}})$ there exists a constant
$C_{p,q}$ that depends on $p$ and $q$ only such that for any centered Gaussian
vector $X$
\[\|X\|_{p}\leq C_{p,q}\|X\|_{q}.\]
\end{Cor}

{\bf Proof.} It is well known that $\|X\|_{p}\leq C_{p}M$, so it is
enough to show that for each $q>r:=\min(-1,-\frac{M^{2}}{4\sigma^{2}})$,
$\|X\|_{q}\geq c_{q}M$. However by (\ref{estdil1}) and Theorem
\ref{estgauss2} we get $\Pr(\|X\|\leq tM)\leq (Ct)^{r}$ and the desired
estimate immediately follows. \hfill $\Box$

\section{Counterexample to Vershynin's conjecture}

Vershynin conjectured that cylinders have the slowest decay of the
Gaussian measure under homothethic shrinking among all centrally symmetric
convex bodies with fixed width and Gaussian measure. Namely, he conjectured
that if $C=B_{2}^{k}(0,w) \times \er^{l}$
(where $B_{2}^{k}(0,w)=\{x\in\er^{k}\colon |x|\leq w\}$)
and $K$ is a centrally symmetric convex
body such that $\gamma(K)=\gamma_{k+l}(C)$ and $w(K)=w=w(C)$ then
$\gamma(tK) \leq \gamma_{k+l}(tC)$ for any $t \in (0,1).$
Certainly, the dimension parameter $l$ is a bit artificial here and one
can easily reduce the problem to the case $l=0.$ This conjecture
seemed naturally related to Conjecture \ref{hipoteza}. Indeed, if a cylinder
$C$ has large width and $\gamma_{k+l}(C)=1/2$ then
$\gamma_{k}(B_{2}^{k}(0,w))=1/2$ and $w$ must be close to $\sqrt{k}.$
On the other hand, $\log_{t} \gamma_{k+l}(tC) \to k$ as $t \to 0.$
Although Conjecture \ref{hipoteza} still seems open, we prove that
Vershynin's cylinder conjecture cannot hold in general.

For simplicity the counterexample will be produced for $k=2$ but one can
easily extend our construction to any $k>2$.

Let us recall that for a set $A$ in $\er^{n}$
\[\gamma_{n}^{+}(A):=\liminf_{t\rightarrow 0}
  \frac{\gamma_{n}(A_{t})-\gamma_{n}(A)}{t}.\]
We begin with two simple and quite standard lemmas.

\begin{lem}
\label{boundary}
There exist positive numbers $w$ and $a$ such that
$\gamma_{2}(B_{2}^{2}(0,w))=\gamma_{1}((-a,a))$ and
$\gamma_{2}^{+}(B_{2}^{2}(0,w)) > \gamma_{1}^{+}((-a,a))$.
\end{lem}

{\bf Proof.}
Let $w$ and $a$ be positive numbers such that
$\gamma_{2}(B_{2}^{2}(0,w))=\gamma_{1}((-a,a)).$ We will prove that
if $a$ is large enough then
$\gamma_{2}^{+}(B_{2}^{2}(0,w)) > \gamma_{1}^{+}((-a,a)).$
Let us recall a standard estimate for the Gaussian tails:
\[
a^{-1}e^{-a^{2}/2}-\int_{a}^{\infty}e^{-s^{2}/2}ds=
\]
\[
\int_{a}^{\infty} \frac{d}{dx} (\int_{x}^{\infty} e^{-s^{2}/2}ds -
x^{-1}e^{-x^{2}/2})dx=
\int_{a}^{\infty}x^{-2}e^{-x^{2}/2}dx,
\]
so that
\[
\int_{a}^{\infty}e^{-s^{2}/2}ds \leq a^{-1}e^{-a^{2}/2}
\]
and
\[
\int_{a}^{\infty}e^{-s^{2}/2}ds \geq a^{-1}e^{-a^{2}/2}-
a^{-3}\int_{a}^{\infty}xe^{-x^{2}/2}dx=(a^{-1}-a^{-3})e^{-a^{2}/2}.
\]
Since
\[
e^{-w^{2}/2}=1-\gamma_{2}(B_{2}^{2}(0,w))=1-\gamma_{1}((-a,a))=
\]
\[
\frac{2}{\sqrt{2\pi}} \int_{a}^{\infty}e^{-s^{2}/2}ds \leq
\frac{2}{\sqrt{2\pi}} a^{-1}e^{-a^{2}/2}
\]
we have $w \geq a+a^{-1}\ln a + o(a^{-1}\ln a)$ as $a \to \infty.$
Therefore for sufficiently large $a$ we have
\[
\gamma_{2}^{+}(B_{2}^{2}(0,w))=we^{-w^{2}/2} >
(a+2a^{-1})e^{-w^{2}/2}=(a+2a^{-1}) \cdot
\frac{2}{\sqrt{2\pi}} \int_{a}^{\infty}e^{-s^{2}/2}ds \geq
\]
\[
\frac{2}{\sqrt{2\pi}} (a+2a^{-1})(a^{-1}-a^{-3})e^{-a^{2}/2} >
\frac{2}{\sqrt{2\pi}} e^{-a^{2}/2}=\gamma_{1}^{+}((-a,a)).
\]

\vspace{-1.5ex}\hfill $\Box$

\begin{lem}
\label{slln}
 If a sequence of positive numbers $u(n)$ satisfies
\[
\liminf_{n \to \infty} \gamma_{n}(B_{2}^{n}(0,u(n))) >0
\]
then $\liminf_{n \to \infty} n^{-1/2}u(n) \geq 1.$ Conversely,
if a sequence of positive numbers $v(n)$ satisfies
$\liminf_{n \to \infty} n^{-1/2}v(n) >1$ then
\[
\lim_{n \to \infty} \gamma_{n}(B_{2}^{n}(0,v(n)))=1.
\]
\end{lem}

{\bf Proof.} Let $g_{1}, g_{2}, \ldots$ be i.i.d. ${\cal N}(0,1)$
random variables. The assertion easily follows from the observation that
\[
\gamma_{n}(B_{2}^{n}(0,r))={\bf P}(\sum_{i=1}^{n}g_{i}^{2} < r^{2})
\]
and the Law of Large Numbers. \hfill $\Box$

\medskip

{\bf Contruction of the counterexample.}
Let $w, a >0$ be such that
$\gamma_{2}(B_{2}^{2}(0,w))=\gamma_{1}((-a,a))$ and
$\gamma_{2}^{+}(B_{2}^{2}(0,w)) > \gamma_{1}^{+}((-a,a))$
and let $t \in (0,1).$ For $x \in (0,w)$ let
$y=y(x)=\frac{x}{2}+\frac{w^{2}}{2x}$ and $s=s(x)=\frac{2w^{2}x}{w^{2}+x^{2}}.$
Note that $x<s<y.$
We define a continuous function $f:[0,\infty) \to [0,w]$ by:\\
$f(r)=w$ if $r \in [0,x],$\\
$f(r)=0$ if $r \geq y$ and\\
$f(r)=\frac{(y-r)w}{y-x}$ if $r \in (x,y).$

Let $K_{n}=K_{n}(x)$ be a {\em flying saucer} body defined by
\[
K_{n}=\{ \xi \in \er^{n}: f((\sum_{i=1}^{n-1}\xi_{i}^{2})^{1/2}) > |\xi_{n}| \}.
\]

Clearly, $K_{n}$ is a convex body contained in the symmetric strip
$S=\{ \xi \in \er^{n}: |\xi_{n}| \leq w \}$. Notice that the line on
the plane $\{(\xi_{1},\xi_{2})\colon \xi_{2}=(y-\xi_{1})w/(y-x)\}$ is
tangent to the ball $B_{2}^{2}(0,w)$ and the tangent point is
$(s,f(s))$. Thus $K_{n}$ contains the inscribed
Euclidean ball $B_{2}^{n}(0,w)$.
Hence the inradius $w(K_{n})=w.$ For sufficiently large $n$ one can choose
$x \in (0,w)$ such that $\gamma_{n}(K_{n})=\gamma_{2}(B_{2}^{2}(0,w)).$
Indeed, for fixed $n,$ the body $K_{n}$ tends to $S$ as $x \to 0^{+}$ and
it tends to $C_{n}=B_{2}^{n-1}(0,w) \times (-w, w)$ as $x \to w^{-}$. Since
$\gamma_{n}(S)>\gamma_{2}(B_{2}^{2}(0,w))$ and
$\gamma_{n}(C_{n}) \to 0$ as $n \to \infty,$ our claim follows by the
continuity. Moreover, $x=x(n) \to 0$ as $n \to \infty.$
Indeed, it suffices to note that $y(n)=y(x(n)) \to \infty$
since
\[
\liminf_{n \to \infty}
\gamma_{n-1}(B_{2}^{n-1}(0,y(n)))\geq
\liminf_{n \to \infty}
\gamma_{n}(B_{2}^{n-1}(0,y(n)) \times (-w,w)) \geq
\]
\[
\gamma_{n}(K_{n})
=\gamma_{2}(B_{2}^{2}(0,w))>0.
\]
From now on we assume that $n$ is large enough and
$x=x(n)$ is such that $\gamma_{n}(K_{n})=\gamma_{2}(B_{2}^{2}(0,w)).$
Let $b, c$ and $d$ be positive numbers such that $a>b>c>d.$
Let $u, v, z \in (x,y)$ be such that $f(u)=b, f(v)=c$ and $f(z)=d.$
Simple calculations show that
\[
u=x+\frac{1}{2}(\frac{w}{x}-\frac{x}{w})(w-b),
\]
\[
v=x+\frac{1}{2}(\frac{w}{x}-\frac{x}{w})(w-c)
\]
and
\[
z=x+\frac{1}{2}(\frac{w}{x}-\frac{x}{w})(w-d).
\]
Hence $v/u \to \frac{w-c}{w-b}>1$ and $v/z \to \frac{w-c}{w-d}<1$
as $n \to \infty.$
Since
\[
K_{n} \subset (\er^{n-1} \times (-b,b)) \cup
(B_{2}^{n-1}(0,u) \times ([b,w) \cup (-w,-b]))
\]
we have
\[
\gamma_{1}((-a,a))=\gamma_{2}(B_{2}^{2}(0,w))=\gamma_{n}(K_{n}) \leq
\gamma_{1}(-b,b)+2\gamma_{1}((b,w))\gamma_{n-1}(B(0,u))
\]
and therefore
\[
\liminf_{n \to \infty} \gamma_{n-1}(B_{2}^{n-1}(0,u)) > 0,
\]
so that $\liminf_{n \to \infty} n^{-1/2}u \geq 1.$
Hence $\liminf_{n \to \infty} n^{-1/2}v \geq \frac{w-c}{w-b}>1$
and consequently
\[
\lim_{n \to \infty} \gamma_{n-1}(B_{2}^{n-1}(0,v))=1.
\]
If $t> \frac{w-c}{w-d}$ then
\[
\liminf_{n \to \infty} \frac{tz}{v} \geq t \frac{w-d}{w-c} >1,
\]
so that $tz \geq v$ for $n$ large enough. Moreover
\[
B_{2}^{n-1}(0,z) \times (-d,d) \subset K_{n},
\]
so
\[
B_{2}^{n-1}(0,tz) \times (-td,td) \subset tK_{n}
\]
and (assuming Vershynin's cylinder conjecture is true) we get that
for any $\frac{w-c}{w-d}<t<1$,
\[
\gamma_{2}(B_{2}^{2}(0,tw))=\gamma_{2}(tB_{2}^{2}(0,w)) \geq
\gamma_{n}(tK_{n}) \geq 
\gamma_{n-1}(B_{2}^{n-1}(0,tz))\gamma_{1}((-td,td)) \geq
\]
\[
\gamma_{n-1}(B_{2}^{n-1}(0,v))\gamma_{1}((-td,td)) \to
\gamma_{1}(-td,td)
\]
as $n \to \infty.$
We have proved $\gamma_{2}(B_{2}(0,tw)) \geq \gamma_{1}(-td,td)$
if only $\frac{w-c}{w-d} <t <1.$ Note however that for a fixed $t \in (0,1)$
one can set $b, c \to a$ and $d \to (a-(1-t)w)/t$ in such way that
$t > \frac{w-c}{w-d}.$ Then $td \to a-(1-t)w$ and we deduce that
\[
\gamma_{2}(B_{2}^{2}(0,tw)) \geq \gamma_{1}((-(a-(1-t)w),a-(1-t)w))
\]
for any $t \in (0,1).$
The above inequality becomes an equality for $t=1$ so by differentiating
at $t=1$ we obtain
\[
\gamma_{2}^{+}(B_{2}^{2}(0,w)) \leq \gamma_{1}^{+}((-a,a)),
\]
contrary to the way in which we chose $w$ and $a$ in the beginning.
This proves that Vershynin's cylinder conjecture cannot be true in general.
  
\vspace{2mm}

{\bf Remark} Note that one really needs some ``extra'' dimensions in the
construction of the counterexample. If we assume that $K \subset \er^{k}$
is a convex symmetric body with inradius $w$ and
$\gamma_{k}(K)=\gamma_{k}(B_{2}^{k}(0,w))$ then obviously
$K$ must be equal to $B_{2}^{k}(0,w)$ up to some boundary points,
so that also $\gamma_{k}(tK)=\gamma_{k}(B_{2}^{k}(0,tw)$ for $t \in (0,1)$
even though the Euclidean ball has {\bf the fastest} decay of the Gaussian
measure under homothetic shrinking among all bodies of the fixed Gaussian
measure.

{\bf Acknowledgments.} Part of this research was done during the
Snowbird AMS-IMS-SIAM Conference "Gaussian Measure and Geometric Convexity"
in July 2004. We would like to thank the organizers for creating excellent
research conditions.

\leftline{Institute of Mathematics, Warsaw University}
\leftline{Banacha 2, 02-097 Warsaw, Poland}
\leftline{e-mail: rlatala@mimuw.edu.pl, koles@mimuw.edu.pl}

\end{document}